\documentclass[11pt,a4paper,reqno]{article} 
\usepackage{amsfonts,amsthm, amscd, epsfig, amsmath, amssymb,enumerate,url,psfrag}

\newtheorem{theorem}{Theorem}
\newtheorem{definition}{Definition}
\newtheorem{proposition}[theorem]{Proposition}

\newenvironment{zeproof}{\vskip 2mm\noindent \textbf{\textit{Proof~}}}
                    {\hfill $\blacksquare$ \vskip 2mm \noindent}

\def\eps{\varepsilon}

\newcommand{\E}{\mathop{\hbox{\rm I\kern-0.17em E}}\nolimits}
\renewcommand{\P}{\mathop{\hbox{\rm I\kern-0.17em P}}\nolimits}
\newcommand{\Zp}{\mathbb{Z}_+}
\newcommand{\Z}{\mathbb{Z}}
\newcommand{\N}{\mathbb{N}}
\newcommand{\T}{\mathcal{T}}
\newcommand{\TT}{\mathbb{T}}

\begin{document} 

\title{Stochastic domination for the last passage percolation model}
  
\author{\textsc{David
    Coupier}\footnote{\texttt{david.coupier@math.univ-lille1.fr}}\\  \textsc{Philippe Heinrich}\footnote{\texttt{philippe.heinrich@math.univ-lille1.fr}}}
\date{Laboratoire Paul Painlev\'e, UMR 8524\\
 UFR de Math\'ematiques, USTL, B\^at. M2\\
59655 Villeneuve d'Ascq Cedex France \\
\texttt{date of submission \today}}

\maketitle

{\bf Abstract:} A competition model on $\mathbb{Z}_+^{2}$ governed by directed last
passage percolation is considered. A stochastic domination argument between
subtrees of the last passage percolation is put forward.  



\section{Introduction}

The directed last passage percolation model goes back to the original work of Rost \cite{Rost} in the case of i.i.d. exponential weights. In this paper, Rost proved a shape theorem for the infected region and exhibited for the first time a link with the one-dimensional totally asymmetric simple exclusion process (TASEP). A background on exlusion processes can be found in the book \cite{Liggett} of Liggett. Since then, this link has been done into details by Ferrari and its coauthors \cite{FGM,FMP,FP} to obtain aymptotic directions and related results for competition interfaces. Other results have been obtained in the case of i.i.d. geometric weights : see Johansson \cite{Joh}. For i.i.d. weights but with general weight distribution, Martin \cite{Martin} proved a shape theorem and described the behavior of the shape function close to the boundary. See also the survey \cite{Martin3}.

\begin{figure}[!ht]
\begin{center}
\psfrag{0}{$(0,0)$}
\includegraphics[width=5.5cm,height=4.5cm]{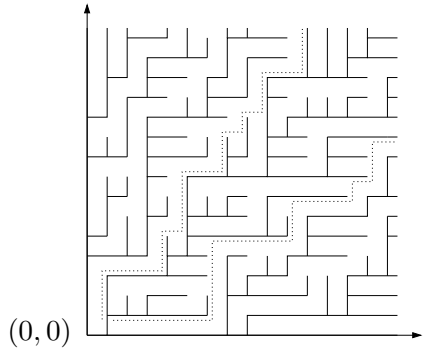}
\end{center}
\caption{\label{fig:lastpassagetree} An example of the last passage percolation tree on the set $[0;15]^{2}$. The subtree rooted at $(1,1)$ is surrounded by dotted lines. Here, the upper dotted line corresponds to the competition interface studied by Ferrari and Pimentel in \cite{FP}.}
\end{figure}

Let us consider $\Omega=[0,\infty)^{\Z^2}$ refered as the configuration space and endowed with a Borel probability measure $\P$. All throughout this paper, $\P$ is assumed translation-invariant : for all $a\in \Z^2$,
$$
\P=\P \circ \tau_a^{-1} ~,
$$
where $\tau_a$ denotes the translation operator on $\Omega$ defined by $\tau_a(\omega)=\omega(a+\cdot)$. This is the only assumption about the probability measure $\P$. We are interested in the behavior of optimal paths from the origin to a site $z\in\Zp^2$. The collection of optimal paths forms the last passage percolation tree $\T$. In this paper, a special attention is paid to the subtree of $\T$ rooted at $(1,1)$: see Figure \ref{fig:lastpassagetree}.

Our goal is to stochastically dominate subtrees of the last passage percolation tree by the one rooted at $(1,1)$. Our results (Theorems \ref{thm2} and \ref{thm3}) essentially rely on elementary properties of the last passage percolation model; its directed nature and the positivity of weights.

The paper is organized as follows. In the rest of this section, optimal paths and the last passage percolation tree $\T$ are precisely defined. The growth property which allows us to compare subtrees of $\T$ is introduced. Theorems \ref{thm2} and \ref{thm3} are stated and commented in Section \ref{theorems}. They are proved in Section \ref{proofs}.

\subsection{Paths, low-optimality, percolation tree}

We will focus on (up-right oriented only) \emph{paths} which can be defined as
sequences (finite or not) $\gamma=(z_0,z_1,\ldots)$ of sites $z_i\in \Z^2$
such that $z_{i+1}-z_i=(1,0)\text{ or }(0,1)$. 

For a given configuration $\omega$, we define the \emph{length} of a path
$\gamma$ as
$$\omega(\gamma)=\sum_{z\in \gamma}\omega(z).$$
If $\Gamma_z$ is the (finite) set of paths from $(0,0)$ to $z$, a path $\gamma\in \Gamma_z$ is  $\omega$-\emph{optimal} if its length $\omega(\gamma)$ is maximal on $\Gamma_z$. The quantity $\max_{\gamma\in\Gamma_z}\omega(\gamma)$ is known as the \emph{last passage time} at $z$. To avoid questions on uniqueness of optimal paths, it is
convenient to call \emph{low-optimal} the optimal path below all the others.

\begin{proposition}
\label{pro1}
Given $\omega\in\Omega$,  each $\Gamma_z$ contains a (unique) low-optimal path denoted by $\gamma_z^\omega$.
\end{proposition}

\begin{zeproof}
 We can assume that $\mathrm{Card}(\Gamma_z)\ge 2$. Given $\omega\in\Omega$, consider two arbitrary
optimal paths $\gamma,\gamma'$ of $\Gamma_z$. If they have no common point (except endpoints $(0,0)$ and $z$), then one path is below
 the other. If $\gamma$ and $\gamma'$  meet in sites, say $u_1,\ldots,u_k$, it's easy to see
 that the path which consists in concatenation of lowest subpaths of
 $\gamma,\gamma'$ between
 consecutive $u_i, u_{i+1}$ is also an optimal path  of $\Gamma_z$. This
 procedure can be (finitely) repeated to reach the low-optimal path of $\Gamma_z$ for the configuration $\omega$.
\end{zeproof}

In literature, optimal paths are generally unique and called \emph{geodesics}. This is the case when $\P$ is a product measure over $\Z^2$ of non-atomic laws. Here, low-optimality ensures uniqueness without particular restriction and Proposition~\ref{pro1} allows then to define the (last passage) percolation tree $\T^\omega$ as the collection of
low-optimal paths $\gamma_z^\omega$  for all
$z\in\Zp^2$. Moreover, the subtree of $\T^\omega$
rooted at $z$ is denoted by $\T_z^\omega$.

\subsection{Growth property}

Let us introduce the set $\TT$ of all substrees of $\T$ : 
$$\TT=\{\T_z^\omega\,:\,z\in \Zp^2,\,\omega\in \Omega\}.$$
For a tree $T\in\TT$, $r(T)$ and $V(T)$ denote respectively its root and its
vertex set.
\begin{definition}
\label{def:growth}
 A subset $A$ of $\TT$ satisfies the \emph{growth
  property} if 
  \begin{equation}
    \label{growth}
    (T\in A , \; T'\in \TT , \; V(T)-r(T)\subset V(T')-r(T')) \; \Longrightarrow \; T'\in A. 
  \end{equation}
\end{definition}
\noindent
For example, if $k\in\Zp\cup\{\infty\}$, the set $\{T\in \TT\,:\,\mathrm{Card}V(T)\ge k\}$ satifies the growth property. But so does not the set 
$$
\{T\in \TT\,:\,T \textrm{ have at least two infinite branches}\} ~.
$$
Indeed, the partial ordering on the set $\TT$ induced by Definition \ref{def:growth} does not take into account the graph structure of trees.

\section{Stochastic domination}
\label{theorems}

The following results compare subtrees of the last passage percolation tree through subsets of $\TT$ satisfying the growth property. 

\begin{theorem}
\label{thm2}
Let $a\in \Zp^{2}$ and a subset $A$ of $\TT$
satisfying the growth property \eqref{growth}. Set also $\Omega^{a}=\left\{\omega\in \Omega:a \textrm{ belongs to
  }\gamma_{a+(1,0)}^\omega\textrm{ and }\gamma_{a+(0,1)}^\omega\right\}$.
 Then, 
$$\P\left(\mathcal{T}_{a+(1,1)}\in A,\Omega^{a}\right)\le
\P\left(\mathcal{T}_{(1,1)}\in A,\tau_a(\Omega^{a})\right).$$ 
  In particular, if $\P$ is in addition a product measure, we have
 $$\P\left(\mathcal{T}_{a+(1,1)}\in A\,\mid\,\Omega^{a}\right)\le
\P\left(\mathcal{T}_{(1,1)}\in A\right),$$
\end{theorem}

To illustrate the meaning of this result, assume that the vertices of $\T_{(1,1)}$ are painted in blue and those of $\T_{(2,0)}$ and $\T_{(0,2)}$ in red. This random coloration leads to a competition of colors. The red area is necessarily unbounded since the model forces every vertice $(x,0)$ or $(0,x)$ with $x\in \{2,3,\ldots\}$ to be red. But the blue area can be bounded. Now, consider  $a\in \Zp^2$ and the same way to color but only in the quadrant $a+\Zp^2$ : this time, the blue area consists of the vertices of $\T_{a+(1,1)}$ and the red one of
$$
\T_{a+(x,0)} \; \mbox{ and } \; \T_{a+(0,x)} , \; \mbox{ for } \; x \geq 2 ~.
$$
Roughly speaking, Theorem \ref{thm2} says that, conditionally to $\Omega^{a}$, the competition is harder for the latter blue area.

The proof of Theorem~\ref{thm2} can be summed up as follows. From a configuration $\omega$, a new one which is a perturbed translation of  $\omega$, namely $\omega^a=\tau_a(\omega)+\eps$ is built in order to satisfy 
$$\T_{a+(1,1)}^\omega=a+\T_{(1,1)}^{\omega^a}.$$
But $\eps$ is chosen such that for $\omega\in \Omega^{a}$, we have $V(\T_{(1,1)}^{\omega^a})\subset V(\T_{(1,1)}^{\tau_a(\omega)})$ and the growth property leads to 
$$\T_{a+(1,1)}^\omega\in A \; \Longrightarrow \; \T_{(1,1)}^{\tau_a(\omega)}\in A.$$
It remains then to use the translation invariance of $\P$ to get the result.

The next result suggests a second stochastic domination argument in the spirit of Theorem~\ref{thm2}.
\begin{theorem}
\label{thm3}
Let $m\in\N$ and a subset $A$ satisfying the growth property \eqref{growth}. Set $\Omega_m=\{\omega\,:\,\gamma_{(m,1)}^\omega=\left((0,0),(1,0),\ldots,(m,0),(m,1)\right)\}$. Then
\begin{equation}
\label{comparaison}
\P(\T_{(m,1)}\in A,\Omega_m)\le \P(\T_{(1,1)}\in A,\Omega_1) ~.
\end{equation}
\end{theorem}

\begin{figure}[!ht]
\begin{center}
\psfrag{1}{\small{$(1,1)$}}
\psfrag{m}{\small{$(m,1)$}}
\includegraphics[width=8cm,height=2.8cm]{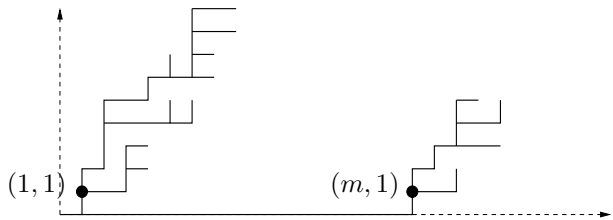}
\end{center}
\caption{\label{fig:comp_arbres} Are represented the subtrees of the last passage percolation tree rooted at sites $(1,1)$ and $(m,1)$, for a configuration $\omega\in\Omega^{1}\cap\Omega^{m}$.}
\end{figure}

Now some comments are needed. 
\begin{itemize}
\item Note that $\Omega_{1}=\{\omega\,:\,\omega(1,0)>\omega(0,1)\}$.
\item It is worth pointing out here that Theorem \ref{thm3} is, up to a certain extend, better than Theorem \ref{thm2}. If $a=(m,0)$ then the events $\Omega^{a}$ and $\Omega_{m}$ are equal and the probablity $\P(\T_{a+(1,1)}\in A,\Omega^{a})$ can be splitted into 
\begin{equation}
\label{bound1thm3}
\P \left( \T_{(m+1,1)} \in A , \; \Omega_{m} , \; \omega(m+1,0) < \omega(m,1) \right)
\end{equation}
and
\begin{equation}
\label{bound2thm3}
\P \left( \T_{(m+1,1)} \in A , \; \Omega_{m} , \; \omega(m+1,0) \geq \omega(m,1) \right) ~.
\end{equation}
On the event $\{\omega(m+1,0)<\omega(m,1)\}$, $\T_{(m+1,1)}$ is as a subtree of $\T_{(m,1)}$. Hence, if $A$ satisfies the growth property (\ref{growth}) then $\T_{(m+1,1)}\in A$ forces $\T_{(m,1)}\in A$. It follows that (\ref{bound1thm3}) is bounded by $\P(\T_{(m,1)}\in A,\Omega_{m})$ which is at most $\P(\T_{(1,1)}\in A,\Omega_{1})$ by Theorem \ref{thm3}. 

On the other hand, $\{\Omega_{m},\omega(m+1,0)\geq\omega(m,1)\}$ is included in $\Omega_{m+1}$. Consequently, \eqref{bound2thm3} is bounded by $\P(\T_{(m+1,1)}\in A,\Omega_{m+1})$, and also by $\P(\T_{(1,1)}\in A,\Omega_{1})$ by Theorem \ref{thm3} again. 

Combining these bounds, we get
$$
\P( \T_{a+(1,1)} \in A , \; \Omega^{a} ) \leq 2 \P( \T_{(1,1)} \in A , \; \Omega_{1} ) ~.
$$
To sum up, whenever $2\P(\T_{(1,1)}\in A,\Omega_{1})$ is smaller than $\P(\T_{(1,1)}\in A)$ (this is the case when $\P$ and $A$ are invariant by the symmetry with respect to the diagonal $x=y$), Theorem \ref{thm2} with $a=(m,0)$ can be obtained as a consequence of Theorem \ref{thm3}.
\item Let us remark that further work seems to lead to the following improvement of Theorem \ref{thm3}: the application
$$
m \; \mapsto \P \left( \T_{(m,1)} \in A \;,\; \Omega^{m} \right)
$$
should be non increasing.
\item Finally, by symmetry, Theorem \ref{thm3} obviously admits an analogous version on the other axis. Roughly speaking, the subtree of the last passage percolation tree rooted at the site $(1,m)$ is stochastically dominated by the one rooted at $(1,1)$.
\end{itemize}

Here are two situations in which Theorem~\ref{thm3}  can be used.\\
An infinite low-optimal path is said non trivial if it does not coincide with one of the two axes $\Zp(1,0)$ and $\Zp(0,1)$. If the set $V(\T_{(1,1)})$ is unbounded (which can be referred as ``coexistence'') then, since each vertex in a subtree has a bounded number of children (in fact, at most 2), the tree $\T_{(1,1)}$ contains an infinite low-optimal path. So, if we set $$Coex=\{\mathrm{Card}V(\T_{(1,1)})=\infty\},$$
then
\begin{equation}
\label{coexistence}
\P \left( Coex\right) > 0 \; \Longrightarrow \; \P \left( \begin{array}{c}
\mbox{there exists a non trivial,} \\
\mbox{infinite low-optimal path}
\end{array} \right) > 0 ~.
\end{equation}
Conversely, assume that $\P(Coex)$ is zero. Since the set
$$\{T\in\TT\,:\, \mathrm{Card}(V(T))=\infty\}$$ 
satisfies the growth property, Theorem \ref{thm3} implies that for all $m \in \N$
$$
\P \left( \mathrm{Card}(V(\T_{(m,1)}))=\infty , \; \Omega_{m} \right) = 0 ~.
$$
Hence, $\P-$a.s., each subtree coming from the axis $\Zp(1,0)$ is finite. This result can be generalized to the two axes $\Zp(1,0)$ and $\Zp(0,1)$ by symmetry. Then, $\P-$a.s., there is no non trivial, infinite low-optimal path and (\ref{coexistence}) becomes an equivalence.

Now, Set 
$$\Delta_n=\{ (x,y) \in \mathbb{Z}_{+}^{2}\,:\, x+y=n \},$$
and let us denote by $\alpha_{n}$ the (random) number of vertices of $\T_{(1,1)}$ meeting $\Delta_n$:
$$
\alpha_{n} = \mathrm{Card} \left( V(\T_{(1,1)}) \cap  \Delta_n\right) ~.
$$
The event $Coex$ can be written $\bigcap_{n\in\N}\{\alpha_{n}>0\}$. We say there is ``strong coexistence'' if
$$
\limsup_{n\to\infty} \frac{\alpha_{n}^{\omega}}{n} > 0 ~.
$$
In a future work, Theorem \ref{thm3} is used so as to give sufficient conditions ensuring strong coexistence with positive probability.

\section{Proofs}
\label{proofs}

\subsection{Proof of Theorem~\ref{thm2}}
Recall that $\gamma_z^\omega$ denotes the low-optimal path from $0$ to $z$ for $\omega$.

\begin{enumerate}[I/]
 \item Let $\omega$ and $\omega+\eps$ be configurations where $\eps$ is a vanishing configuration except on
the axes $\Zp (1,0)$ and $\Zp (0,1)$ i.e $\eps(x,y)=0$ if $xy\ne 0$.
We shall show that if $\eps$ satifies
\begin{equation}
 \label{condeps}
\eps(0,1)+\eps(0,2)\ge \eps (1,0) \textrm{ and } \eps(1,0)+\eps(2,0)\ge \eps (0,1),
\end{equation}
then 
\begin{equation}
  \label{eq:prec}
 V(\T_{(1,1)}^{\omega+\eps})\subset V(\T_{(1,1)}^\omega). 
\end{equation}
Let $z\in V(\T_{(1,1)}^{\omega+\eps})$. By definition, $z$ is a vertex  such that the low-optimal path $\gamma_z^{\omega+\eps}$
contains $(1,1)$ ; we have to show that $z\in V(\T_{(1,1)}^{\omega})$ i.e. the low-optimal path $\gamma_z^\omega$
also contains $(1,1)$. 

There is nothing to prove if  $\gamma_z^{\omega+\eps}= \gamma_z^\omega$, so we assume that $\gamma_z^{\omega+\eps}\ne
\gamma_z^\omega$.  By additivity of
$\omega\mapsto\omega(\gamma)$ and optimality of $\gamma_z^{\omega+\eps}$ and $\gamma_z^\omega$, we have
\begin{eqnarray}
 \eps(\gamma_z^\omega) &= &
(\omega+\eps)(\gamma_z^\omega)-\omega(\gamma_z^\omega)\nonumber\\
&\le & (\omega+\eps)(\gamma_z^{\omega+\eps})-\omega(\gamma_z^\omega)\nonumber\\
&=&  \omega(\gamma_z^{\omega+\eps})+ \eps(\gamma_z^{\omega+\eps})-\omega(\gamma_z^\omega)\nonumber\\
&\le & \eps(\gamma_z^{\omega+\eps}).\label{strict0}
\end{eqnarray}
Note also that by low-optimality of  $\gamma_z^{\omega+\eps}$ and $\gamma_z^\omega$, we have 
$$\left((\omega+\eps)(\gamma_z^{\omega+\eps})=(\omega+\eps)(\gamma_z^\omega)\textrm{ and } \omega(\gamma_z^{\omega+\eps})=\omega(\gamma_z^{\omega}) \right)\Longrightarrow \gamma_z^{\omega+\eps}= \gamma_z^\omega.$$
Since $\gamma_z^{\omega+\eps}$ and $\gamma_z^\omega$ are different, this allows us to strengthen \eqref{strict0} : 
\begin{equation}
  \label{eq:strict}
\eps(\gamma_z^\omega)<  \eps(\gamma_z^{\omega+\eps}).
\end{equation}

Now, it's got to be one thing or the other:
\begin{itemize}
 \item Either the site $(1,0)$ belongs to $\gamma_z^{\omega+\eps}$. In this case, the right hand side of  \eqref{eq:strict} which becomes  $\eps(0,0)+\eps (1,0)$ and  the strict inequality imply that  
$\gamma_z^\omega$ can not run through $(1,0)$. Moreover, if $\gamma_z^\omega$ ran through $(0,2)$, we would have
$$\eps(0,0)+\eps (0,1)+\eps (0,2)\le \eps(\gamma_z^\omega)<\eps(0,0)+\eps (1,0),$$
but this would be in contradiction with \eqref{condeps}. We conclude that $\gamma_z^\omega$ must run through $(0,1)$ and $(1,1)$. 
\item  Or the site $(0,1)$ belongs to $\gamma_z^{\omega+\eps}$, and
  symmetrically we conclude that, if \eqref{condeps} hold, $\gamma_z^\omega$ runs through $(1,0)$ and $(1,1)$.
\end{itemize}
To sum up, if $\eps$ satisfies \eqref{condeps} then \eqref{eq:prec} holds. The conditions \eqref{condeps} can be understood as follows; the second one (for example) prevents the set $V(\T_{(2,0)})$ to drop vertices in favour of $V(\T_{(1,1)})$ passing from $\omega$ to $\omega+\eps$.
\item For given configuration $\omega\in\Omega$ and site $a\in\Zp^2$, we construct a new
  configuration $\omega^{a}$ such that 
  \begin{equation}
    \label{eq:cong}
\forall z\in \Zp^2,\quad
\omega^{a}\left(\gamma_z^{\omega^{a}}\right)=
\omega\left(\gamma_{a+z}^\omega\right).
  \end{equation}
The idea of the constuction is to translate $\omega$ from $a$ to the origin and to modify then weights on the axes : more precisely, set
$$\omega^{a}(z)=
\begin{cases}
 \omega(\gamma_{a}^\omega) & \text{if $z=(0,0)$;}\\
\omega(\gamma_{a+z}^\omega)- \omega(\gamma_{a+z-(1,0)}^\omega)&
\text{if $z=(x,0) $ with $x\in \N$;}\\
\omega(\gamma_{a+z}^\omega)- \omega(\gamma_{a+z-(0,1)}^\omega)&
\text{if $z=(0,y)$ with $y \in \N$;}\\
\omega(a+z) & \text{otherwise.}\\
\end{cases}
$$
Let $\bar{z}$ be the latest site of $a+(\Zp(1,0)\cup\Zp(0,1))$ whereby $\gamma_{a\!+\!z}^{\omega}$ passes. The configuration $\omega^{a}$ is defined so as to the last passage time to $\bar{z}$ for $\omega$ is equal to the last passage time to $\bar{z}-a$ for $\omega^{a}$, i.e. $\omega(\gamma_{\bar{z}}^{\omega})=\omega^{a}(\gamma_{\bar{z}-a}^{\omega^{a}})$. Combining with $\omega^{a}(\cdot)=\omega(a+\cdot)$ on $a+\mathbb{N}^{2}$, the identity (\ref{eq:cong}) follows. Finally, by low-optimality, the translated path $a+\gamma_{z}^{\omega^{a}}$ coincides with the restriction of $\gamma_{a+z}^\omega$ to the quadrant $a+\mathbb{N}^{2}$. See Figure \ref{fig:omegaz0}.

\begin{figure}[!ht]
\begin{center}
\psfrag{a}{\small{$(0,0)$}}
\psfrag{b}{\small{$a$}}
\psfrag{c}{\small{$\bar{z}$}}
\psfrag{d}{\small{$a\!+\!z$}}
\psfrag{e}{\small{$\bar{z}\!-\!a$}}
\psfrag{f}{\small{$z$}}
\includegraphics[width=11cm,height=4.2cm]{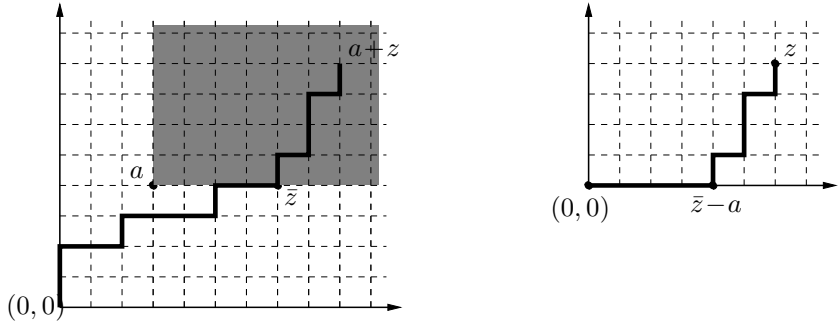}
\end{center}
\caption{\label{fig:omegaz0} To the left, the low-optimal path to $a+z$ for a given configuration $\omega$ is represented. Let us denote by $\bar{z}$ the latest site of $a+(\Zp(1,0)\cup\Zp(0,1))$ whereby $\gamma_{a\!+\!z}^{\omega}$ passes. To the right, the low-optimal path to $z$ for the corresponding configuration $\omega^{a}$ is represented.}
\end{figure}

In particular, we can write with
some abuse of notation 
\begin{equation}
  \label{eq:tree}
a+\mathcal{T}_{(1,1)}^{\omega^{a}}=\mathcal{T}_{a+(1,1)}^\omega.
\end{equation}
The induction formula's 
\begin{equation}
 \label{induc}
\omega(\gamma_u^\omega)=\max(\omega(\gamma_{u-(1,0)}^\omega),\omega(\gamma_{u-(0,1)}^\omega))+\omega(u),
\end{equation} 
allows to rewrite the configuration $\omega^{a}$:
\begin{equation}
\label{eps}
 \omega^{a}=\tau_a(\omega)+\eps,
\end{equation}
where $\eps$ is defined on the axes by
\begin{eqnarray}
\eps(0,0) &= & \max\left(\omega(\gamma_{a-(1,0)}^\omega),\omega(\gamma_{a-(0,1)}^\omega)\right) \nonumber\\
\eps(x,0) & = & \max\left(0,\omega(\gamma_{a+(x,-1)}^\omega)-\omega(\gamma_{a+(x-1,0)}^\omega)\right)\quad (x\in\N) \label{ex}\\
\eps(0,y) & = & \max\left(0,\omega(\gamma_{a+(-1,y)}^\omega)-\omega(\gamma_{a+(0,y-1)}^\omega)\right)\quad (y\in\N) \label{ey}\\
\eps(x,y) &=& 0 \quad \textrm{otherwise} ~. \nonumber
\end{eqnarray}

\item Consider $a$ and 
$\Omega^{a}$ as in the statement of Theorem~\ref{thm2}. Let $\omega\in
\Omega^{a}$ so that the length $\omega(\gamma_{a}^\omega)$ is bigger than
$\omega(\gamma_{a+(1,-1)}^\omega)$ and
$\omega(\gamma_{a+(-1,1)}^\omega)$. It follows from \eqref{ex} and \eqref{ey}
that $\eps(1,0)=\eps(0,1)=0$. Conditions \eqref{condeps} are then trivially satisfied so that
\eqref{eq:prec} holds for $\omega$ and also for $\tau_a(\omega)$:
$$V(\T_{(1,1)}^{\tau_a(\omega)+\eps})\subset V(\T_{(1,1)}^{\tau_a(\omega)})~.$$
Combined with \eqref{eq:tree} and \eqref{eps}, this leads to
$$V(\mathcal{T}_{a+(1,1)}^\omega)-a\subset V(\T_{(1,1)}^{\tau_a(\omega)}).~$$
Now, if $A$ satisfies the growth property \eqref{growth} then
$$\mathcal{T}_{a+(1,1)}^\omega\in A \; \Longrightarrow \; \T_{(1,1)}^{\tau_a(\omega)}\in A.$$
To summarize, we have $\{\mathcal{T}_{a+(1,1)}\in A\}\subset \tau_{a}^{-1}\{\mathcal{T}_{(1,1)}\in A\}$ on  $\Omega^{a}$,
and since $\P$ is translation-invariant and $\Omega^a=\tau_a^{-1}(\tau_a(\Omega^{a}))$, we conclude that
\begin{eqnarray*}
 \P( \mathcal{T}_{a+(1,1)}\in A,\Omega^{a}) &\le&  \P\left(\tau_{a}^{-1}\{\mathcal{T}_{(1,1)}\in A\},\Omega^{a}\right) \\
&= &\P(\mathcal{T}_{(1,1)}\in A,\tau_a(\Omega^{a})).
\end{eqnarray*}
The first part of Theorem \ref{thm2} is proved. In order to prove the second part, let us assume $\P$ is a product measure. It suffices to remark the events $\tau_a(\Omega^{a})$ which means both low-optimal paths from $-a$ to $(1,0)$ and $(0,1)$ run through the origin, and $\{\mathcal{T}_{(1,1)}\in A\}$ are independent. Actually, the random variable $\omega(0,0)$ is the only weight of $\mathbb{Z}_+^{2}$ of which $\tau_a(\Omega^{a})$ depends on, and it is involved in all optimal paths coming from the origin. So, it does not affect the event $\{\mathcal{T}_{(1,1)}\in A\}$.
\end{enumerate}

\subsection{Proof of Theorem~\ref{thm3}}

\begin{enumerate}[I/]
\item Let $\omega\in\Omega_{1}$ and $\omega+\eps$ be two configurations where $\eps$ is a vanishing configuration except on the axis $\Zp(0,1)$: $\eps(x,y)=0$ whenever $x>0$. We also assume that $\omega$ and $\eps$ verify $\omega(1,0)>\omega(0,1)+\eps(0,1)$ (i.e. $\omega+\eps\in\Omega_{1}$). The goal of the first step consists in stating:
\begin{equation}
\label{eq:prec2}
V(\T_{(1,1)}^{\omega+\eps}) \subset V(\T_{(1,1)}^\omega ) ~. 
\end{equation}
Let $z$ be a vertex such that the low-optimal path $\gamma_z^{\omega+\eps}$ contains $(1,1)$. If the low-optimal paths $\gamma_z^\omega$ and $\gamma_z^{\omega+\eps}$ are different then it follows as for \eqref{eq:strict}):
$$
\eps(\gamma_z^\omega) < \eps(\gamma_z^{\omega+\eps}) ~.
$$
Henceforth, the condition $\omega+\eps\in\Omega_{1}$ implies that $\gamma_z^{\omega+\eps}$ runs through $(1,0)$ and leads to a contradiction:
$$
\eps(0,0) \leq \eps(\gamma_z^\omega) < \eps(\gamma_z^{\omega+\eps}) = \eps(0,0) ~.
$$
So, $\gamma_z^\omega$ and $\gamma_z^{\omega+\eps}$ are equal, which implies $z$ is a vertex of $\T_{(1,1)}^\omega$. Relation (\ref{eq:prec2}) is proved. It is worth to note that condition $\omega+\eps\in\Omega_{1}$ ensures that the random interface between sets $V(\T_{(1,1)})$ and $V(\T_{(2,0)})$ remains unchanged if we add $\eps$ to $\omega$. Hence, the set $V(\T_{(1,1)})$ can only decrease.

\item Let $\omega$ be a configuration and $b=(m-1,0)$. In the spirit of the proof of Theorem \ref{thm2}, a configuration $\omega^{b}$ is built by translating $\omega$ by vector $-b$ and preserving the last passage percolation tree structure. The right construction is the following:
$$
\omega^{b}(z) =\begin{cases}
\omega(\gamma_{b}^\omega) & \text{if $z=(0,0)$;}\\
\omega(\gamma_{b+z}^\omega)- \omega(\gamma_{b+z-(0,1)}^\omega)&
\text{if $z\in \{0\}\times \N$;}\\
\omega(b+z) & \text{otherwise}.
\end{cases}
$$
By construction, the configuration $\omega^{b}$ satisfies $\omega^{b}(\gamma_z^{\omega^{b}})=\omega(\gamma_{b+z}^\omega)$, for all $z\in\mathbb{Z}_{+}^{2}$. Thus, we can deduce from low-optimality:
\begin{equation}
\label{eq:tree2}
b + \mathcal{T}_{(1,1)}^{\omega^{b}} = \mathcal{T}_{b+(1,1)}^{\omega} ~.
\end{equation}
The induction formula's \eqref{induc} allows to write for all $z \in \mathbb{Z}_{+}^{2}$,
$$
\omega^{b}(z) = \tau_b(\omega)(z) + \eps(z) ~,
$$
with 
$$
\eps(z) =\begin{cases}
\omega(\gamma_{b-(1,0)}^\omega) & \text{if $z=(0,0)$;}\\
\max \left( \omega(\gamma_{b+z-(1,0)}^\omega)-\omega(\gamma_{b+z-(0,1)}^\omega) , 0 \right)&
\text{if $z\in \{0\}\times \N$;}\\
 0 & \text{otherwise.}
\end{cases}
$$
Besides,
\begin{eqnarray}
\label{equiv:m1}
\omega \in \Omega_{m} & \iff & \omega(\gamma_{b+(0,1)}^{\omega}) < \omega(\gamma_{b+(1,0)}^{\omega}) \nonumber\\
& \iff & \omega^{b}(0,1) + \omega(\gamma_{b}^{\omega}) < \omega(\gamma_{b}^{\omega}) + \omega(b+(1,0)) \nonumber\\
& \iff & \omega^{b}(0,1) < \omega^{b}(1,0) \nonumber\\
& \iff & \omega^{b} \in \Omega_{1} ~.
\end{eqnarray} 

\item Given $\omega\in\Omega_{m}$, equivalence \eqref{equiv:m1} implies $\omega^{b}=\tau_b(\omega)+\eps\in\Omega_{1}$. As a by-product, we have $\tau_b(\omega)\in\Omega_{1}$ and from (\ref{eq:prec2}) and (\ref{eq:tree2}), we deduce
$$
V(\mathcal{T}_{b+(1,1)}^{\omega})-b \subset V(\T_{(1,1)}^{\omega^{b}}) \subset V(\T_{(1,1)}^{\tau_b(\omega)}) ~.
$$
If $A\subset\TT$ satisfies the growth property (\ref{growth}) then
$$
\left(\omega \in \Omega_{m} \; \mbox{ and } \; \mathcal{T}_{(m,1)}^{\omega} \in A \right)\; \Longrightarrow \; \left(\tau_b(\omega) \in \Omega_{1} \; \mbox{ and } \; \mathcal{T}_{(1,1)}^{\tau_b(\omega)} \in A \right) ~.
$$
Finally, \eqref{comparaison} easily follows from the translation invariance of the probability measure $\P$.
\end{enumerate}

\bibliographystyle{plain}
\bibliography{stochdom_bibli}

\end{document}